\documentclass[12pt]{amsart}
\usepackage{amssymb}
\usepackage{graphics}
\usepackage{colordvi}

\newtheorem{theorem}{Theorem}

\newtheorem{defn}[theorem]{Definition}
\newtheorem{ex}[theorem]{Example}

\newenvironment{definition}{\begin{defn}\rm}{\end{defn}}
\newenvironment{example}{\begin{ex}\rm}{\end{ex}}
\newcounter{FNC}[page]
\def\newfootnote#1{{\addtocounter{FNC}{1}$^\fnsymbol{FNC}$%
     \let\thefootnote\relax\footnotetext{$^\fnsymbol{FNC}$#1}}}
\newcommand{\C}{\mathbb{C}}
\newcommand{\Z}{\mathbb{Z}}

\newcommand{\Ln}{\raisebox{-2pt}{$\includegraphics{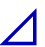}_n$}}

\title[Horn recursion for Schur $P$- and $Q$- functions]{The Horn recursion for
  Schur $P$- and $Q$- functions:\\ Extended Abstract}  

\author{Kevin Purbhoo}
\address{Department of Mathematics\\
         The University of British Columbia\\
         1984 Mathematics Road\\
         Vancouver, B.C., V6T 1Z2\\
         CANADA}
\email{kevinp@math.ubc.ca}
\urladdr{http://www.math.ubc.ca/~kevinp/}
\author{Frank Sottile}
\address{Department of Mathematics\\
         Texas A\&M University\\
         College Station\\
         TX \ 77843\\
         USA}
\email{sottile@math.tamu.edu}
\urladdr{http://www.math.tamu.edu/\~{}sottile}
\thanks{Work of Sottile supported by the Clay Mathematical Institute and NSF
         CAREER grant DMS-0134860} 
\thanks{Work of Purbhoo supported by NSERC} 



\begin{document}

\begin{abstract}
 A consequence of work of Klyachko and of Knutson-Tao is the Horn recursion to
 determine when a Littlewood-Richardson coefficient is non-zero.
 Briefly, a Littlewood-Richardson coefficient is non-zero if and only if it
 satisfies a collection of Horn inequalities which are indexed by smaller
 non-zero Littlewood-Richardson coefficients.
 There are similar Littlewood-Richardson numbers for Schur $P$- and $Q$-
 functions. 
 Using a mixture of combinatorics of root systems, combinatorial linear algebra
 in Lie algebras, and the geometry of certain cominuscule flag
 varieties, we give Horn recursions to determine when these
 other Littlewood-Richardson numbers are non-zero.
 Our inequalities come from the usual Littlewood-Richardson numbers,
 and while we give two very different Horn recursions, they have the same 
 sets of solutions.
 Another combinatorial by-product of this work is a new Horn-type recursion for
 the usual Littlewood-Richardson coefficients.

\end{abstract}


\maketitle

\section*{Introduction}
The Littlewood-Richardson numbers $a^\lambda_{\mu,\nu}$ for partitions
$\lambda,\mu,\nu$  are important in many areas of mathematics.
For example, they are the structure constants of several related rings with
distinguished bases:
the ring of symmetric functions with its basis of Schur
functions, the representation ring of $\mathfrak{sl}_n$ with its basis of
irreducible highest weight modules, the external representation ring of the
tower of symmetric groups with its basis of irreducible modules, and
the cohomology ring of the Grassmannian with its basis of Schubert
classes~\cite{Fu97,Mac,St99}.
The combinatorics of Littlewood-Richardson numbers are extremely interesting
and now we have many formulas for them, including the original
Littlewood-Richardson formula~\cite{LR1934}. 
Despite this deep and prolonged interest in Littlewood-Richardson numbers, one
of the most fundamental questions about them was not asked until about a decade 
ago: 
\[
  \mbox{When is $a^\lambda_{\mu,\nu}$ non-zero?}
\]

This question came from (of all places) a problem in linear algebra concerning
the possible eigenvalues of a sum of hermitian matrices.
The answer to this problem is given by the Horn inequalities:
the eigenvalues which can and do occur are the solutions to a set of linear
inequalities, and the inequalities themselves come from non-negative integral
eigenvalues solving this problem for smaller matrices.

The same inequalities answer our question about Littlewood-Richardson numbers.
A Littlewood-Richard\-son number $a^\lambda_{\mu,\nu}$ is non-zero if and only if
the triple of partitions $(\lambda,\mu,\nu)$ satisfy certain linear inequalities,
and the inequalities themselves come from triples indexing  smaller non-zero
Littlewood-Richardson coefficients.
This description is a consequence of work of Klyachko~\cite{Klyachko}
which linked eigenvalues of sums of hermitian matrices, highest weight modules
of $\mathfrak{sl}_n$, and the Schubert calculus for Grassmannians, and then
Knutson and Tao's proof~\cite{KT99} of Zelevinsky's Saturation
Conjecture~\cite{Ze99}.
This work implies Horn's Conjecture~\cite{Ho62}
about the eigenvalues of sums of Hermitian matrices.
These results have wide implications in mathematics (see the
surveys~\cite{Fu98,Fu00}) and have raised many new and evocative
questions. 
For example, the Horn inequalities give the answer to questions in 
several different realms of mathematics (representation theory, 
combinatorics, Schubert calculus, eigenvalues), 
but it was initially mysterious why any one of these questions should
have a recursive answer,
as the proofs travelled through so many other parts of mathematics.

Another question, which was the {\it point de d\'epart} for the results we
describe here, is the following:
are there related numbers whose non-vanishing has a similar recursive answer? 
Our main result is a recursive answer for the non-vanishing of the
analogs of Littlewood-Richardson coefficients for Schur $P$-functions, and the
same for Schur $Q$-functions.
We give one set of inequalities which determine non-vanishing for the
$P$-functions and a different set of inequalities for the $Q$-functions. 
Because each Schur $P$-function is a non-zero multiple of a corresponding Schur
$Q$-function, a Littlewood-Richardson number for $P$-functions is non-zero
if and only if the corresponding number for $Q$-functions is non-zero, and
thus our two sets of inequalities have the  same sets of solutions.

Another consequence of our work is a new set of recursive Horn-type 
inequalities for the ordinary Littlewood-Richardson numbers
$a^\lambda_{\mu,\nu}$. 
While these new inequalities are clearly related to the ordinary Horn
inequalities, they are definitely quite different.
(We explain this below.)

Before we define some of these objects and give the different
recursions, we remark that our  results were proved using a mixture
of the combinatorics of root systems, combinatorial linear algebra in Lie algebras,
and the geometry of certain \Blue{{\it cominuscule}} flag varieties $G/P$.
Cominuscule flag varieties are (almost all of) the flag varieties whose
geometrically defined Bruhat order (which is the Bruhat order on the cosets
$W/W_P$ of the Weyl groups) is a distributive lattice.

The alert reader will notice that these inequalities for Schur $P$- and
$Q$-functions are not strictly recursive because they are indexed by 
ordinary Littlewood-Richardson numbers which are non-zero.
The reason for the term recursive is that the inequalities
stem from a geometric recursion among all cominuscule flag varieties which
is not evident from viewing only the subclass corresponding to, say the Schur
$Q$-functions.

This abstract does not dwell on the geometry, but rather on the fascinating
combinatorial consequences of these recursions. 
The last section of this extended abstract does however give a broad view of
some of the key geometric ideas which underly our recursion.
The results here are proved in the forthcoming preprint by the authors, 
``The recursive nature of the cominuscule Schubert calculus''.

\section{The classical Horn recursion}
For more details and definitions concerning the various flavors of Schur
functions that arise here, we recommend the book of
Macdonald~\cite{Mac}.
\Blue{{\it Schur functions $S_\lambda$}} are symmetric functions indexed by 
\Blue{{\it partitions $\lambda$}}, which are weakly decreasing sequences of
nonnegative integers,
$\lambda\colon\lambda_1\geq\lambda_2\geq\dotsb\geq\lambda_n\geq 0$.
The Schur function $S_\lambda$ is homogeneous of degree 
$|\lambda|:=\lambda_1+\dotsb+\lambda_n$.
Schur functions form a basis for the $\Z$-algebra of symmetric functions.
Thus there are integral \Blue{{\it Littlewood-Richardson}} numbers
$a^\lambda_{\mu,\nu}$ defined by the identity
\[
  S_\mu\cdot S_\nu\ =\ \sum_\lambda \Blue{a^\lambda_{\mu,\nu}}\,S_\lambda\,.
\]
Homogeneity gives the necessary relation $|\lambda|=|\mu|+|\nu|$ for 
$a^\lambda_{\mu,\nu}\neq 0$.

A partition $\lambda$ is represented by its \Blue{{\it diagram}}, which is a
left-justified array of boxes in the positive quadrant with $\lambda_i$ boxes in
row $i$. 
Thus
\[
   (4,2,1) \ \leftrightarrow\  \raisebox{-9pt}{\includegraphics{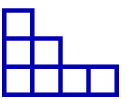}}
\]
Partitions are partially ordered by the inclusion of their diagrams.
Let $n\times m$ be the rectangular partition with $n$ parts, each of size
$m$. 

For $\lambda\subset n\times m$ ($\lambda_1\leq m$ and $\lambda_{n+1}=0$), define
$\lambda^c$ to be the partition obtained from the set-theoretic difference 
$n\times m - \lambda$ of diagrams (placing $\lambda$ in the upper right corner
of $n\times m$).
Thus we have
\[
  \begin{picture}(62,42)
  \put(0,0){\includegraphics{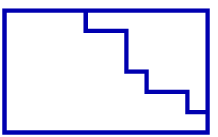}}
  \put(10,11){$\lambda^c$}  \put(46,23){$\lambda$}
  \put(25,41){$m$}  \put(64,19){$n$}
   \put(78,20){.}
  \end{picture}
\]

The Horn-type inequalities we give are naturally stated in terms of 
\Blue{{\it symmetric Littlewood-Richardson numbers}}.
For $\lambda,\mu,\nu\subset n\times m$, define\
 \begin{eqnarray*}
  \Blue{a_{\lambda,\mu,\nu}} &:=&
    \mbox{Coefficient of $S_{n\times m}$ in $S_\lambda S_\mu S_\nu$}\\
    &\ = &
     \mbox{Coefficient of $S_{\lambda^c}$ in $S_\mu S_\nu$}
     \ =\ a^{\lambda^c}_{\mu,\nu}\,.
 \end{eqnarray*}
We say that a triple of partitions $\lambda,\mu,\nu\subset n\times m$
is \Blue{{\it feasible}} if $a_{\lambda,\mu,\nu}\neq 0$.
This convenient terminology comes from geometry.
 
\begin{definition}
 Suppose that $\lambda\subset n\times m$ and $\alpha\subset r\times(n-r)$,
 where  $0< r < n$.
 Define 
\[
   I_n(\alpha)\ :=\ \{n-r+1-\alpha_1,\ n-r+2-\alpha_2,\ \dotsc,\ n-\alpha_r\}\,.
\]
 Draw $\lambda$ in the upper right corner of the $n\times m$ rectangle, and
 number the rows Cartesian-style.
 Define $\Blue{|\lambda|_\alpha}$ to be the number of boxes that remain in
 $\lambda$ after crossing out the rows indexed by $I_n(\alpha)$.
\end{definition}

\begin{example}\label{Ex:ClassicHorn}
 Suppose that $n=7$, $m=8$, and $r=3$, 
 and we have $\lambda=8654310$ and $\alpha=311$.
 Then the set $I_7(\alpha)$ is 
\[
   \{7-3+1-3,\ 7-3+2-1,\ 7-3+3-1\}\ =\ \{2, 5, 6\}\,.
\]
 If we place $\lambda$ in the upper-right corner of the rectangle $7\times
 8$ and cross out the rows indexed by $I_7(\alpha)$,
\[
 \begin{picture}(137,95)(0,-2)
   \put(2,  3){1}   \put(2,17){\Blue{2}}   \put(2, 31){3}
   \put(2, 45){4}   \put(2,59){\Blue{5}}   \put(2, 73){\Blue{6}}
   \put(2, 87){7} 

   \put(12,-1){\includegraphics{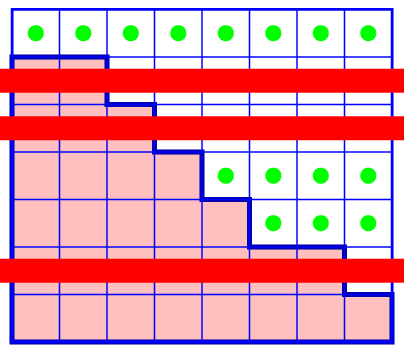}}
  \end{picture}
\]
 we count the dots \Green{{$\bullet$}} to see that  $|\lambda|_\alpha=15$.
\end{example}

\begin{theorem}[Classical Horn Recursion: Klyachko~\cite{Klyachko}, Knutson-Tao~\cite{KT99}]
\ 

 A triple $\lambda,\mu,\nu\subset n\times m$ is feasible if and only if\/
  $|\lambda| + |\mu| + |\nu| = nm$, and 
\[
  |\lambda|_\alpha + |\mu|_\beta + |\nu|_\gamma\ \leq \ 
   (n-r)m
\]
for all feasible triples $\alpha,\beta,\gamma\subset r\times(n-r)$ and for 
all\/ $0< r< n$.
\end{theorem}

The first condition, $|\lambda| + |\mu| + |\nu| = nm$, is due to homogeneity.

\section{Symmetric Horn recursion}

Since replacing a partition $\lambda$ by its conjugate $\lambda^t$
(interchanging rows with columns) induces an involution on the algebra of
symmetric functions, there is a version of the Horn recursion where one crosses
out columns instead of rows.
It turns out that there are necessary inequalities obtained by crossing
out both rows and columns, including possibly a different number of each.
The cominuscule recursion reveals a sufficient subset of these.

\begin{definition}\label{Definition:NewHorn}
 Let $0< r < \min\{n,m\}$ and 
 suppose that $\lambda\subset n\times m$, $\alpha\subset r\times(n-r)$, and 
 we have another partition $\alpha'\subset r\times(m-r)$.
 Define $I_n(\alpha)$ as before, and set
\[
   I_m(\alpha')\ :=\ \{ \Blue{m}-r+1-\alpha'_1,\ 
       \Blue{m}-r+2-\alpha'_2,\ \dotsc,\ \Blue{m}-\alpha'_r\}\,.
\]
 Draw $\lambda$ in the upper right corner of the $n\times m$ rectangle and cross
 out the rows indexed by $I_n(\alpha)$ and the columns indexed by $I_m(\alpha')$.
 Define $\Blue{|\lambda|_{\alpha,\alpha'}}$ to be the number of boxes that
 remain in $\lambda$.
\end{definition}

\begin{example}
 We use the same data as in Example~\ref{Ex:ClassicHorn},  and set
 $\alpha'=410$.
 Then 
\[
   I_8(\alpha')\ =\ \{8-3+1-4,\ 8-3+2-1,\ 8-3+3-0\}\ =\ \{2, 6, 8\}\,.
\]
 If we now cross out the rows indexed by $I_7(\alpha)$ and the columns indexed
 by $I_8(\alpha')$, 
\[
  \begin{picture}(125,120)(0,-2)
   \put( 10,0){1}   \put( 24,0){\Blue{2}}   \put( 38,0){3}
   \put( 52,0){4}   \put( 66,0){5}          \put( 80,0){\Blue{6}}
   \put( 94,0){7}   \put(108,0){\Blue{8}}

   \put(4,11){\includegraphics{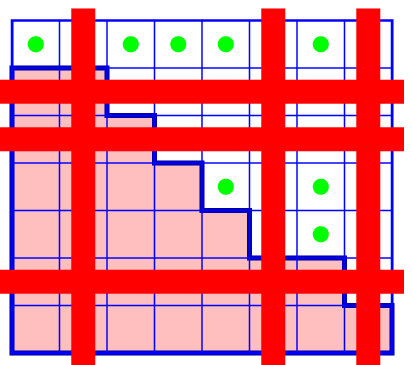}}
 \end{picture}
\]
 we count the dots \Green{{$\bullet$}} to see that 
 $|\lambda|_{\alpha,\alpha'}=8$. 
\end{example}

\begin{theorem}[Symmetric Horn Recursion]
\ 

 A triple $\lambda,\mu,\nu\subset n\times m$ is feasible if and only if\/
  $|\lambda| + |\mu| + |\nu| = nm$, and 
\[
  |\lambda|_{\alpha,\alpha'} + |\mu|_{\beta,\beta'} + |\nu|_{\gamma,\gamma'}
   \ \leq \    (m-r)(n-r)
\]
 for all pairs of feasible triples 
 $\alpha,\beta,\gamma\subset r\times(n-r)$ and
 $\alpha',\beta',\gamma'\subset r\times(m-r)$ and for 
 all\/ $0< r< \min\{m,n\}$.
\end{theorem}

\section{Schur $P$- and $Q$- functions}

The algebra of symmetric functions has a natural odd subalgebra which
comes from its structure as a combinatorial Hopf algebra~\cite{CHA}.
This algebra was first studied by Schur in the context of the theory of
projective representations of the symmetric group.
This odd subalgebra has a pair of distinguished bases, the 
\Blue{{\it Schur $P$-functions}} and  the 
\Blue{{\it Schur $Q$-functions}}.
They are indexed by strict partitions, which are strictly decreasing sequences
of positive integers 
$\lambda\colon \lambda_1>\lambda_2>\dotsb>\lambda_k>0$.
They are proportional: $Q_\lambda=2^k P_\lambda$, where $\lambda$ has $k$ parts.

We have Littlewood-Richardson coefficients $c^\lambda_{\mu,\nu}$ and 
$d^\lambda_{\mu,\nu}$ indexed by triples of strict partitions and defined by the
identities
\[
  Q_\mu\cdot Q_\nu\ =\ \sum_\lambda \Blue{c^\lambda_{\mu,\nu}}\,Q_\lambda
  \qquad\mbox{ and }\qquad
  P_\mu\cdot P_\nu\ =\ \sum_\lambda \Blue{d^\lambda_{\mu,\nu}}\,P_\lambda\ .
\]
Combinatorial formulas for these numbers were given in work of Worley~\cite{Wo},
Sagan~\cite{Sa}, and Stembridge~\cite{Stem89}.

Let $\Ln\colon n>n{-}1>\dotsb>2>1$ be the strict partition of staircase shape.
Then $\lambda\subset\Ln$ if $\lambda_1\leq n$.
If $\lambda\subset\Ln$, define
$\lambda^c$ to be the partition obtained from the set-theoretic difference 
$\Ln - \lambda$ of diagrams (placing $\lambda$ in the upper right corner
of $\Ln$).
Thus we have
\[
   \begin{picture}(70,65)
  \put(0,0){\includegraphics{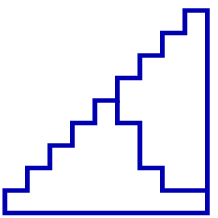}}
    \put(22,10){$\lambda^c$}  \put(45,30){$\lambda$}
    \put(18,37){$n$}
   \end{picture}
\]

As before, the Horn-type inequalities are naturally stated in terms of 
symmetric Littlewood-Richardson numbers.
For $\lambda,\mu,\nu\subset \Ln$, define
\begin{eqnarray*}
\Blue{c_{\lambda,\mu,\nu}} &:=&
  \mbox{Coefficient of $Q\!_{\includegraphics{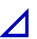}_n}$ in 
             $Q_\lambda Q_\mu Q_\nu$}\\
  &\ = &
   \mbox{Coefficient of $Q_{\lambda^c}$ in $Q_\mu Q_\nu$}
     \ =\ c^{\lambda^c}_{\mu,\nu}\,.
\end{eqnarray*}
A triple of strict partitions $\lambda,\mu,\nu\subset n\times m$
is \Blue{{\it feasible}} if $c_{\lambda,\mu,\nu}\neq 0$.

We similarly define symmetric Littlewood-Richardson numbers
$\Blue{d_{\lambda,\mu,\nu}}$ for the Schur $P$-functions.
Since the two bases are proportional, the corresponding coefficients are as
well.
In particular the sets of triples $\lambda,\mu,\nu$ for which the corresponding
coefficients are feasible are the same.
Nevertheless, we give two very different sets of inequalities which
determine the feasibility of these numbers, arising from the 
different geometric origins of Schur $Q$-functions and Schur $P$-functions.

\begin{definition}\label{Def:typeC}
 Let $0< r < n$ and 
 suppose that $\lambda\subset \Ln$ is a strict partition and $\alpha\subset
 r\times(n-r)$ is an ordinary partition. 
 Draw $\lambda$ in the upper right corner of the staircase $\Ln$.
 Number the inner corners $1,2,\dotsc,n$ and, for 
 each number in $I_n(\alpha)$, cross out the row and column emanating from that
 inner corner.
 Then let  $\Blue{[\lambda]_\alpha}$ be the number of boxes that remain in
 $\lambda$.
\end{definition}

 \begin{definition}\label{Def:typeD}
 Let $0< r < \Blue{n{+}1}$ and 
 suppose that $\lambda\subset \Ln$ is a strict partition and 
 $\alpha\subset r\times(\Blue{n{+}1}{-}r)$ is an ordinary partition. 
 Draw $\lambda$ in the upper right corner of the staircase $\Ln$.
 Number the outer corners $1,2,\dotsc,n,\Blue{n{+}1}$ and for 
 each number in $I_{\Blue{n+1}}(\alpha)$, cross out the row and column emanating
 from the corresponding outer corner.
 Then let  $\Blue{\{\lambda\}_\alpha}$ be the number of boxes that remain in
 $\lambda$. 
\end{definition}

\begin{example}\label{Ex:CD-Horn}
 For example, suppose that $n=8$ and $r=4$, 
 we have $\lambda=8643$ and $\alpha=4220$.
 Then 
 \begin{eqnarray*}
   I_8(\alpha)&=&
    \{8-4+1-4,\ 8-4+2-2,\ 8-4+3-2,\ 8-4+4-0\}\ =\ \{1, 4, 5, 8\}\,,\\
   I_9(\alpha)&=&
    \{9-4+1-4,\ 9-4+2-2,\ 9-4+3-2,\ 9-4+4-0\}\ =\ \{2, 5, 6, 9\}\,.
 \end{eqnarray*}
 and if we place $\lambda$ in the upper-right corner of the rectangle 
 \raisebox{-2pt}{$\includegraphics{Ln.eps}_8$} and cross out the rows and
 columns emanating from the inner corners indexed by $I_8(\alpha)$,
 we see that $[\lambda]_\alpha=6$.
 If we instead cross out the rows and
 columns emanating from the outer corners indexed by $I_9(\alpha)$,
 we see that $\{\lambda\}_\alpha=5$.
 The two diagrams are shown in Figure~\ref{Fig:CD}, on the left and right,
 respectively. 
\begin{figure}[htb]
\[
  \begin{picture}(117,115)(-4,15)
   \put(-4, 25){\Blue{1}} \put( 9, 38){2} \put( 22, 51){3}
   \put(35, 64){\Blue{4}} \put( 48, 77){\Blue{5}} \put(61,90){6}
   \put(74,103){7} \put(87,118){\Blue{8}}  
   \put(-4,105){$[\lambda]_\alpha=6$}
   \put( 4,10){\includegraphics{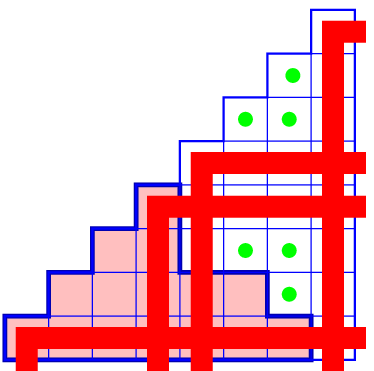}}
  \end{picture}
  \qquad\qquad\qquad
  \begin{picture}(120,120)(0,15)
   \put(-4, 20){1} \put(  9,33){\Blue{2}} \put( 22, 46){3}
   \put(35, 59){4} \put( 48,72){\Blue{5}} \put(61,85){\Blue{6}}
   \put(74, 98){7} \put(87,111){8} \put(100,124){\Blue{9}}
   \put(-4,105){$\{\lambda\}_\alpha=5$}
   \put(12,10){\includegraphics{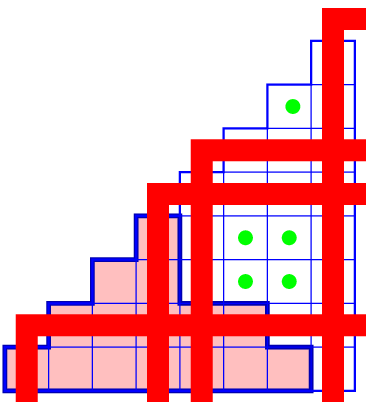}}
  \end{picture}
\]
\caption{ Computation of $[\lambda]_\alpha=6$ and of $\{\lambda\}_\alpha=5$}
\label{Fig:CD}
\end{figure}
\end{example}

Note that the homogeneity of the multiplication of Schur $P$-functions and Schur
$Q$-functions implies that 
 \begin{equation}\label{Eq:homog}
   |\lambda| + |\mu| + |\nu|\ =\ \left|\Ln\right|\ =\ \binom{n+1}{2}\,,
 \end{equation}
 is necessary for a triple $\lambda,\mu,\nu\subset\Ln$ to be feasible.

\begin{theorem}[Horn recursion for Schur $P$- and $Q$-functions]
\ 

 A triple $\lambda,\mu,\nu\subset\Ln$ of strict partitions is feasible
 if and only if
 one of the following two equivalent conditions hold:

\noindent {\rm (1)} 
    The homogeneity condition~\eqref{Eq:homog} holds, and for all feasible
    $\alpha,\beta,\gamma\subset r\times(n-r)$ and  all $0<r<n$, we have
\[
    [\lambda]_\alpha + [\mu]_\beta + [\nu]_\gamma\ \leq\ \binom{n+1-r}{2}\,,
\]
\ \ \  or else 

\noindent {\rm(2)} 
    The homogeneity condition~\eqref{Eq:homog} holds, and for all feasible
    $\alpha,\beta,\gamma\subset r\times(\Blue{n+1}-r)$ and 
       all $0<r<\Blue{n+1}$ with $r$ even, we have
\[
   \{\lambda\}_\alpha + \{\mu\}_\beta + \{\nu\}_\gamma\ \leq\ \binom{n+1-r}{2}\,.
\]
\end{theorem}

\section{Remarks on the geometry of the proof}

We first give some general idea of the ingredients in our proof, and then
explain a little bit of the relation of this geometry to the combinatorics
given here.

A flag manifold $G/P$ ($G$ is a reductive algebraic group and $P$ is a parabolic
subgroup) has a Bruhat decomposition into Schubert cells indexed by
cosets $W/W_P$, where $W$ is the Weyl group of $G$ and $W_P$ that of $P$.
The closures of the Schubert cells are the 
Schubert varieties, and cohomology classes associated to them 
(Schubert classes) form bases for the cohomology ring of $G/P$.
Standard results in geometry show that the structure constants (generalized
Littlewood-Richardson numbers) are the number of points in triple intersections
of general translates of Schubert varieties (and hence are non-negative). 

If a structure constant is non-zero, then any triple
intersection of corresponding Schubert varieties (not just a general
intersection) is non-empty, and general intersections are transverse.
Conversely, if a structure constant is zero, then any corresponding
general intersection is empty, and a non-empty intersection is 
never transverse.
The key idea is to replace the difficult question
on non-emptiness of a general intersection of Schubert varieties by the easier
question of the transversality of a (not completely general) intersection.
This was used by one of us to show transversality of intersections in the
Grassmannian of lines~\cite{So97a}, but its use to study the Horn problem is due to
Belkale~\cite{Be02}, who first gave a geometric proof of the Horn recursion for the
Grassmannian. 

This idea transfers the analysis from the flag manifold $G/P$ to its tangent space 
$T_pG/P$ at a given point $p$.
In fact, all of our diagrams are just pictures of the root-space decompositions
of $T_pG/P$ for the corresponding varieties.
In our proof, we consider three Schubert varieties which
contain the point $p$, and then move them independently by
the stabilizer $P$ of $p$ so that they are otherwise general.
If it is possible to move the three tangent spaces inside of $T_pG/P$ so that
they meet transversally, then the triple is feasible, and if not, then it is not.

This explains where cominuscule flag manifolds come in.
The maximal reductive, or Levi, subgroup $L$ of the 
parabolic group $P$ acts on the tangent space $T_pG/P$ to $G/P$ at that
point. Our arguments (moving the tangent spaces to Schubert varieties around by
elements of $L$) require that $L$ act on $T_pG/P$ with finitely many orbits, and
this is one characterization of cominuscule flag manifolds $G/P$.

It also explains why there is a recursion.
The argument requires us to consider the stabilizer $Q$ in $L$ of a certain
linear subspace of $T_pG/P$---the tangent space to an orbit of $L$ through a
general point of the intersection of general translates of the tangents to the
three Schubert varieties.
Then the Schubert calculus inside of the smaller flag manifold $L/Q$ is used to
analyze the transversality of that triple intersection.
Fortunately, the flag manifold $L/Q$ is itself cominuscule, which is the source
of our recursion.

We briefly illustrate these comments on the three flag manifolds that arose in
this extended abstract.

\subsection{The classical Grassmannian}
Let $Gr(n,m{+}n)$ be the Grassmannian of $n$-planes in $m+n$ space.
The general linear group $GL(m+n,\C)$ acts on $Gr(n,m{+}n)$.
If $H\in Gr(n,m{+}n)$ then $T_HGr(n,m{+}n)$ may be identified with the set of $n$
by $m$ matrices (more precisely with $\mbox{Hom}(H,\C^{m+n}/H)$).
The Levi subgroup is the group $GL(n,\C)\times GL(m,\C)$ which acts linearly on
the rows and columns of $n$ by $m$ matrices.
The orbits of this group are simply matrices of a fixed rank, $r$, and the 
subgroup $Q$ is the stabilizer of a pair $(K,K')$, where 
$K\subset H$ and $K'\subset \C^{m+n}/H$ both have dimension $r$.
This explains why in Definition~\ref{Definition:NewHorn}, the number of rows crossed
out equals the number of columns crossed out.
The smaller cominuscule flag variety $L/Q$ is the product of two Grassmannians,
$Gr(r,n)\times Gr(r,m)$.

The Schubert varieties of $Gr(n,m{+}n)$ are indexed by partitions $\lambda$ which
fit in the $n\times m$ rectangle, and its cohomology ring is the algebra of
Schur functions with these restricted indices.

\subsection{The Lagrangian Grassmannian}
 Fix a non-degenerate alternating bilinear (symplectic) form on
 $\C^{2n}$.
 Let $LG(n)$ be the set of maximal isotropic (Lagrangian) subspaces in $\C^{2n}$,
 each of which has dimension $n$.
 This is the quotient of the symplectic group by the parabolic subgroup
 corresponding to the long root, $LG(n)=Sp(2n,\C)/P_0$.

 Since $H\in LG(n)$ is isotropic the symplectic form identifies $\C^{2n}/H$
 with the dual of $H$, and $T_HLG(n)$ is identified with the space of quadratic
 forms on $H$.
 The Levi subgroup is the general linear group $GL(H)$.
 Identifying $H$ with $\C^n$, the Levi becomes $GL(n,\C)$ and 
 $T_HLG(n)$ is the set of $n\times n$ symmetric matrices.
 (Symmetric matrices are parametrized by their weakly lower triangular parts,
 which correspond to the staircase shape $\Ln$ 
 where the order of the columns has been reversed.)
 The general linear group acts simultaneously on the rows and columns of
 symmetric matrices.
 The orbits are simply symmetric matrices of a fixed rank, $r$, and the 
 subgroup $Q$ is the stabilizer of the null space of such a matrix.
 The smaller cominuscule flag variety $L/Q$ is the Grassmannian
 $G(r,n)$.

The Schubert varieties of $LG(n)$ are indexed by strict partitions $\lambda$ which
fit inside the staircase $\Ln$, and its cohomology ring is the algebra of
Schur $Q$-functions with these restricted indices.

\subsection{The Orthogonal Grassmannian}
 Fix a non-degenerate symmetric bilinear form on $\C^{2n+2}$.
 The set of maximal isotropic subspaces (each of which has dimension $n+1$) 
 of $\C^{2n+2}$ has two isomorphic components.
 Let $OG(n+1)$ be one of these components.
 This is the quotient of the even orthogonal group by a parabolic subgroup $P$
 corresponding to one of the roots at the fork in the Dynkin diagram,
 $OG(n+1)=SO(2n+2,\C)/P$.

 If $H\in OG(n+1)$ is isotropic, then $T_HOG(n+1)$ is identified with the space of
 alternating forms on $H$.
 The Levi subgroup is the general linear group $GL(H)$.
 Identifying $H$ with $\C^{n+1}$, then the Levi becomes $GL(n+1,\C)$ and 
 $T_HOG(n+1)$ is the set of $(n+1)\times(n+1)$ anti-symmetric matrices.
 (Anti-symmetric matrices are parametrized by their lower triangular parts,
 and these strictly lower triangular matrices correspond to the
 staircase shape where the order of the columns has been reversed.)
 The general linear group acts simultaneously on the rows and columns of
 anti-symmetric matrices.
 The orbits are simply anti-symmetric matrices of a fixed rank.
 However, and this comes from the details of the proof and the roots of
 $SO(2n+2,\C)$, the subgroup $Q$ is the stabilizer of an even-dimensional 
 subspace of $H$.
 The smaller cominuscule flag variety $L/Q$ is the Grassmannian
 $G(r,n+1)$, where $r$ is even.

The Schubert varieties of $OG(n{+}1)$ are indexed by strict partitions $\lambda$ which
fit inside the staircase $\Ln$, and its cohomology ring is the algebra of
Schur $P$-functions with these restricted indices.


\providecommand{\bysame}{\leavevmode\hbox to3em{\hrulefill}\thinspace}
\providecommand{\MR}{\relax\ifhmode\unskip\space\fi MR }
\providecommand{\MRhref}[2]{%
  \href{http://www.ams.org/mathscinet-getitem?mr=#1}{#2}
}
\providecommand{\href}[2]{#2}

 
\end{document}